\documentclass[11pt,reqno]{amsart}
\usepackage{graphicx}
\usepackage{verbatim}
\usepackage{textcomp}
\usepackage{amssymb}
\usepackage{cite}
\usepackage{amsmath}
\usepackage{latexsym}
\usepackage{amscd}
\usepackage{amsthm}
\usepackage{mathrsfs}
\usepackage{xypic}
\usepackage{bm}
\usepackage{url}
\usepackage{hyperref}

\newtheorem{thm}{Theorem}[section]

\newtheorem{lem}[thm]{Lemma}

\theoremstyle{definition}

\theoremstyle{remark}
\newtheorem{rem}{Remark}[section]
\numberwithin{equation}{section}
\begin{document}
\title[Some rigidity characterizations on critical metrics]
{Some rigidity characterizations on critical metrics for quadratic curvature
functionals}

\author{Guangyue Huang}
\address{Department of Mathematics, Henan Normal
University, Xinxiang 453007, P.R. China} \email{hgy@henannu.edu.cn }

\thanks{The research of author is supported by NSFC (Nos. 11371018, 11671121).}

\maketitle

\begin{abstract}

We study closed $n$-dimensional manifolds of which the metrics
are critical for quadratic curvature functionals involving the Ricci curvature, the scalar curvature and the Riemannian curvature tensor on the space of Riemannian metrics with unit volume. Under some additional integral conditions, we classify such manifolds. Moreover, under some curvature conditions, the result that a critical metric must be Einstein is proved.

\end{abstract}



\section{Introduction}

Let $\mathscr{M}_1(M^{n})$ be the space of equivalence classes of
smooth Riemannian metrics of volume one on closed Riemannian
manifold $M^{n}$, $n\geq3$. A well-known example of a Riemannian functional is
the Einstein-Hilbert functional
\begin{equation*}
\mathcal{H}=\int_{M}R\,dv
\end{equation*}
on $\mathscr{M}_1(M^{n})$, where $R$ denotes the scalar curvature. By a direct calculation, it is easy to see that Einstein metrics are critical for the functional $\mathcal{H}$.
In this paper, we are interested in studying the functional
\begin{equation}\label{1-Sec-1}
\mathcal{F}_{t,s}(g)\,=\,\int_M |Ric|^2\,dv+t\int_M
R^2\,dv+s\int_M |{\rm Rm}|^2\,dv,
\end{equation}
where $t,s$ are real constants, $Ric$ and $Rm$ denote the Ricci curvature and the Riemannian curvature tensor, respectively. It is easy to observe from \eqref{2-Sec-2} that every Einstein metric is critical for $\mathcal{F}_{t,0}$. In \cite{Catino2015}, Catino considered the curvature functional $\mathcal{F}_{t,0}$ and obtained some conditions on the geometry of $M^n$ such that critical metrics of $\mathcal{F}_{t,0}$ are Einstein. Certainly, there exist critical metrics which are not necessarily Einstein (for
instance, see \cite[Chapter 4]{Bess2008} and \cite{Lamontagne1998}). For some
development in this direction, see
\cite{HL2004,Gursky2015,Gursky2001,Lamontagne1994,Anderson1997,Hu2010}
and the references therein.

The authors in \cite{Barros2017} show that locally conformally flat critical metrics for $\mathcal{F}_{t,s}$ with $n+4(n-1)t+4s=0$ and some additional conditions are space form metrics (see \cite[Theorem 3]{Barros2017}). In this paper, we will give some new characterizations on critical metrics for $\mathcal{F}_{t,s}$ on $\mathscr{M}_1(M^{n})$ with $n+4(n-1)t+4s\neq0$ without assumption that $M^n$ is locally conformally flat. In order to state our results, throughout this paper, we denote by $E$ the traceless Ricci tensor. Our first result reads as follows:

\begin{thm}\label{thm1-1}
Let $M^{n}$ be a closed manifold of dimension $n\neq 4$ with
positive scalar curvature and $g$ be a critical metric for
$\mathcal{F}_{t,s}$ on $\mathscr{M}_1(M^{n})$ with $s\neq-\frac{n-2}{4}$.

(1) If $n+4(n-1)t+4s>0$, then we have
\begin{equation}\label{1For-th-1}\aligned
\int_M&\Big[\frac{1+4s}{2}|C|^2+\frac{(n-4)|n-2+4s|}{(n-2)\sqrt{n(n-1)}}|E|^3+\frac{4-3n-8s-2n(n-1)t}{n(n-1)}R |E|^2\\
&+\frac{|n-2+8s|}{\sqrt{2(n-2)(n-1)}}|W||E|^2+2|s||W|^2|E|\Big]\,dv\geq0
\endaligned\end{equation}
with equality occurs if and only if $M^{n}$ is either Einstein, or
$M^{n}$ is isometrically covered by
$\mathbb{R}\times \mathbb{S}^{n-1}$ with a product metric. In the latter case, we have
$R=\sqrt{n(n-1)}|E|$ and $\frac{n-4}{n-2}|n-2+4s|+4-3n-8s-2n(n-1)t=0$.

(2) If $n+4(n-1)t+4s<0$, then we have
\begin{equation}\label{1For-th-2}\aligned
\int_M&\Big[\frac{1+4s}{2}|C|^2-\frac{(n-4)|n-2+4s|}{(n-2)\sqrt{n(n-1)}}|E|^3+\frac{4-3n-8s-2n(n-1)t}{n(n-1)}R |E|^2\\
&-\frac{|n-2+8s|}{\sqrt{2(n-2)(n-1)}}|W||E|^2-2|s||W|^2|E|\Big]\,dv\leq0
\endaligned\end{equation}
with equality occurs if and only if $M^{n}$ is either Einstein, or
$M^{n}$ is isometrically covered by
$\mathbb{R}\times \mathbb{S}^{n-1}$ with a product metric. In the latter case, we have
$R=\sqrt{n(n-1)}|E|$ and $-\frac{n-4}{n-2}|n-2+4s|+4-3n-8s-2n(n-1)t=0$.

\end{thm}

Under the condition that $E$ is a Codazzi tensor, we also obtain the following integral inequality:

\begin{thm}\label{thm1-2}
Let $M^{n}$ be a  closed manifold of dimension $n\neq 4$ with
positive scalar curvature and $g$ be a critical metric for
$\mathcal{F}_{t,s}$ on $\mathscr{M}_1(M^{n})$ with $s\neq-\frac{n-2}{4}$.
If $1+4s\geq0$ and $E$ is a Codazzi tensor, then we have
\begin{equation}\label{3For-th-1}\aligned
\int_{M}\Big[&-\frac{|2(n-2)+4ns|}{\sqrt{2(n-2)(n-1)}}|W||E|^{\frac{n-2}{n}}-2|s||W|^2|E|^{-\frac{2}{n}}\\
&-\frac{|4s(n^2-3n+4)+4(n-2)|}{(n-2)\sqrt{n(n-1)}}|E|^{\frac{2(n-1)}{n}}\\
&+\frac{4-2n-2n(n-1)t+4(n-2)s}{n(n-1)}R|E|^{\frac{n-2}{n}}\Big]\leq0,
\endaligned\end{equation} and equality occurs if and only if $M^{n}$ is either Einstein, or $M^{n}$ is isometrically covered by $\mathbb{R}^1\times
\mathbb{S}^{n-1}$ with a product metric. In the latter case, we have
$-|4s(n^2-3n+4)+4(n-2)|+(n-2)[4-2n-2n(n-1)t+4(n-2)s]=0$ and $R=\sqrt{n(n-1)}|E|$.

\end{thm}

Next, we provide a rigidity result on critical metrics for $\mathcal{F}_{t,s}$ on $\mathscr{M}_1(M^{n})$ under the assumption that $M^{n}$ is locally conformally flat.

\begin{thm}\label{thm1-3}
Let $M^{n}$ be a locally conformally flat closed manifold of dimension $n\neq4$ with
positive scalar curvature and $g$ be a critical metric for
$\mathcal{F}_{t,s}$ on $\mathscr{M}_1(M^{n})$.

(1) For $n=3$, $t,s$ satisfy
\begin{equation}\label{3lem-Form-11}
\frac{24(1+4s)^2}{(1+6t-2s)^2}|E|^2<R^2
\end{equation}
and
\begin{equation}\label{3lem-Form-12}
\begin{cases}1+2t+2s>0\\
1+4s>0\\
3+8t+4s>0\\
1+6t-2s<0\\
\frac{24(1+4s)^2}{(1+6t-2s)^2}>\frac{(5+16t+4s)^2}{2(3+8t+4s)(1+2t+2s)},
\end{cases}
\end{equation}
then $M^3$ is of constant positive sectional curvature.

(2) For $n\geq5$, $t,s$ satisfy
\begin{equation}\label{3lem-Form-21}
\frac{n(n-1)[4s(n^2-3n+4)+4(n-2)]^2}{(n-2)^2[4-2n-2n(n-1)t+4(n-2)s]^2}|E|^2<R^2
\end{equation}
and
\begin{equation}\label{3lem-Form-22}
\begin{cases}
1+2t+2s>0\\
(n-2)+4s>0\\
n+4(n-1)t+4s<0\\
4-2n-2n(n-1)t+4(n-2)s>0\\
-\Big(\frac{1+4s}{2}+\frac{(n-1)(n-4)(n-2+4s)(1+2t+2s)}{n(n-2)[n+4(n-1)t+4s]}\Big)^2
\frac{2n^2[n+4(n-1)t+4s]}{(n-2)(n-4)(1+4s)(1+2t+2s)[(n-2)+4s]}\\
\leq\frac{n(n-1)[4s(n^2-3n+4)+4(n-2)]^2}{(n-2)^2[4-2n-2n(n-1)t+4(n-2)s]^2},
\end{cases}
\end{equation}

then $M^n$ is of constant positive sectional curvature.

\end{thm}

\begin{rem}

When $s=0$, our Theorem \ref{thm1-1} can provide a concise version.

\end{rem}

\begin{rem}

If the Riemannian curvature tensor of $M^n$ is harmonic (that is, $R_{ijkl,l}=0$), then from the second Bianchi identity, we have
\begin{equation}\label{1-Rem-1}
R_{ki,j}-R_{kj,i}=R_{ijkl,l}=0,
\end{equation}
which shows that the Ricci curvature is Codazzi. Moreover, we have
\begin{equation}\label{1-Rem-2}
R_{,i}=R_{kk,i}=R_{ik,k}=\frac{1}{2}R_{,i},\end{equation}
which shows that $R_{,i}=0$ and the scalar curvature is constant.
Therefore, for $M^n$ with harmonic curvature tensor, we have that $E$ is a Codazzi tensor.

\end{rem}

\begin{rem}

For $n=3$, the existence of $t,s$ satisfying \eqref{3lem-Form-12} is clear. For example, when $t=0$ and $s\geq\frac{17}{25}$. On the other hand, when $s=0$, Catino has proved (see \cite[Theorem 1.5]{Catino2015}) that for $n=3$, if $t\in
[-\frac{1}{3},-\frac{1}{6})$ and
$$|E|<-\frac{1+6t}{2\sqrt{6}}R,$$
then $M^3$ is of
constant positive sectional curvature. Hence, our case of $n=3$ in Theorem \ref{thm1-3} generalizes those
of Catino in \cite{Catino2015}. For $n\geq5$ and $s=0$, the authors in \cite{HC2017} have shown that there exist $\varepsilon_n>0$
and $\eta_n>0$ such that for $-\frac{1}{2}<t<-\varepsilon_n$ or $-\eta_n<t<-\frac{n}{4(n-1)}$,
the inequalities in \eqref{3lem-Form-22} are true.

\end{rem}

\section{Preliminaries}

For $n\geq3$, the Weyl curvature tensor and the Cotton tensor are defined by
\begin{equation}\label{2-Sec-4}\aligned
W_{ijkl}=&R_{ijkl}-\frac{1}{n-2}(R_{ik}g_{jl}-R_{il}g_{jk}
+R_{jl}g_{ik}-R_{jk}g_{il})\\
&+\frac{R}{(n-1)(n-2)}(g_{ik}g_{jl}-g_{il}g_{jk})\\
=&R_{ijkl}-\frac{1}{n-2}(E_{ik}g_{jl}-E_{il}g_{jk}
+E_{jl}g_{ik}-E_{jk}g_{il})\\
&-\frac{R}{n(n-1)}(g_{ik}g_{jl}-g_{il}g_{jk}),
\endaligned\end{equation}
and
\begin{equation}\label{2-Sec-5}\aligned
C_{ijk}=&R_{kj,i}-R_{ki,j}-\frac{1}{2(n-1)}(R_{,i}g_{jk}-R_{,j}g_{ik})\\
=&E_{kj,i}-E_{ki,j}+\frac{n-2}{2n(n-1)}(R_{,i}g_{jk}-R_{,j}g_{ik}),
\endaligned\end{equation}
respectively. From the definition of the Cotton tensor, it is easy to see
\begin{align}\label{2-Sec-6}
C_{ijk}=-C_{jik},\qquad g^{ij}C_{ijk}=g^{ik}C_{ijk}=g^{jk}C_{ijk}=0
\end{align}
and
\begin{align}\label{2-Sec-7}
C_{ijk,k}=0,\qquad C_{ijk}+C_{jki}+C_{kij}=0.
\end{align}
For $n\geq4$, the divergence of the Weyl curvature tensor is related to the Cotton
tensor by
\begin{equation}\label{2-Sec-8}
-\frac{n-3}{n-2}C_{ijk}=W_{ijkl,l}.
\end{equation}

It is well known that $W_{ijkl}=0$ holds naturally on $(M^3,g)$, and $(M^3,g)$ is
locally conformally flat if and only if $C_{ijk}=0$. For $n\geq4$,
$(M^n,g)$ is locally conformally flat if and only if $W_{ijkl}=0$.

First, we recall the following result proved by Catino in \cite{Catino2015} (see \cite[Proposition 6.1]{Catino2015}):

\begin{lem}\label{2lem-1}
Let $M^n$ be a closed manifold of dimension $n\geq3$. A metric $g$ is critical for
$\mathcal{F}_{t,s}$ on $\mathscr{M}_1(M^{n})$ if and only if it satisfies the equations
\begin{equation}\label{2-Sec-2}\aligned
(1+4s)\Delta E_{ij}=&(1+2t+2s)R_{,ij}-\frac{1+2t+2s}{n}(\Delta R)g_{ij}-2(1+2s)R_{ikjl}E_{kl}\\
&-\frac{2+2nt-4s}{n}RE_{ij}+\frac{2}{n}(|E|^2+s|{\rm Rm}|^2)g_{ij}\\
&-2sR_{ikpq}R_{jkpq}+4sE_{ik}E_{jk}
\endaligned\end{equation}
and
\begin{equation}\label{2-Sec-3}\aligned
{[n+4(n-1)t+4s]}\Delta R=&(n-4)(|Ric|^2+tR^2+s|{\rm Rm}|^2-\lambda)\\
=&(n-4)\Big[s|W|^2+\frac{n-2+4s}{n-2}|E|^2\\
&+\frac{n-1+n(n-1)t+2s}{n(n-1)}R^2-\lambda\Big],
\endaligned\end{equation}
where $\lambda=\mathcal{F}_{t,s}(g)$.

\end{lem}

From the Lemma \ref{2lem-1}, we can obtain the following

\begin{lem}\label{2lem-2}
Let $M^n$ be a closed manifold of dimension $n\geq3$. If the metric $g$ is critical for
$\mathcal{F}_{t,s}$ on $\mathscr{M}_1(M^{n})$, then
\begin{equation}\label{2lem-form1}\aligned
\frac{1+4s}{2}\Delta |E|^2=&(1+4s)|\nabla E|^2+(1+2t+2s)R_{,ij}E_{ij}\\
&-\frac{2(n-2)+4ns}{n-2}W_{ikjl}E_{kl}E_{ij}-2sW_{ikpq}W_{jkpq}E_{ij}\\
&+\frac{4s(n^2-3n+4)+4(n-2)}{(n-2)^2}E_{ik}E_{kj}E_{ji}\\
&+\frac{4-2n-2n(n-1)t+4(n-2)s}{n(n-1)}R|E|^2
\endaligned\end{equation}
and
\begin{equation}\label{2lem-form2}\aligned
(1+4s)C_{kij,k}=&\frac{n+4(n-1)t+4s}{2(n-1)}\mathring{R}_{,ij}-\frac{n-2+8s}{n-2}W_{ikjl}E_{kl}\\
&-2s\Big(W_{ikpq}W_{jkpq}-\frac{1}{n}|W|^2g_{ij}\Big)\\
&-\frac{(n-4)(n-2+4s)}{(n-2)^2}\Big(E_{ik}E_{jk}-\frac{1}{n}|E|^2g_{ij}\Big)\\
&+\frac{4-3n-8s-2n(n-1)t}{n(n-1)}R E_{ij},
\endaligned\end{equation}
where $\mathring{R}_{,ij}=R_{,ij}-\frac{R}{n}g_{ij}$.

\end{lem}

\proof
Using the formula \eqref{2-Sec-4}, we can derive
\begin{equation}\label{2-Sec-15}\aligned
E_{kl}R_{ikjl}=&E_{kl}W_{ikjl}+\frac{1}{n-2}(|E|^2g_{ij}-2E_{ik}E_{jk})-\frac{1}{n(n-1)}R E_{ij},
\endaligned\end{equation}
\begin{equation}\label{2-Sec-16}\aligned
R_{ikpq}R_{jkpq}=&W_{ikpq}W_{jkpq}+\frac{4}{n-2}W_{ikjl}E_{kl}+\frac{2(n-4)}{(n-2)^2}E_{ik}E_{jk}\\
&+\frac{2}{(n-2)^2}|E|^2g_{ij}+\frac{2}{n^2(n-1)}R^2g_{ij}+\frac{4}{n(n-1)}RE_{ij}
\endaligned\end{equation}
and
\begin{equation}\label{2-Sec-17}\aligned
|{\rm Rm}|^2=&|W|^2+\frac{4}{n-2}|E|^2+\frac{2}{n(n-1)}R^2.
\endaligned\end{equation}
Therefore, \eqref{2-Sec-2} can be written as
\begin{equation}\label{add2-Sec-2}\aligned
(1+4s)\Delta E_{ij}=&(1+2t+2s)\mathring{R}_{,ij}-2(1+2s)R_{ikjl}E_{kl}\\
&-\frac{2+2nt-4s}{n}RE_{ij}+\frac{2}{n}(|E|^2+s|{\rm Rm}|^2)g_{ij}\\
&-2sR_{ikpq}R_{jkpq}+4sE_{ik}E_{jk}\\
=&(1+2t+2s)\mathring{R}_{,ij}-\frac{2(n-2)+4ns}{n-2}W_{ikjl}E_{kl}-2sW_{ikpq}W_{jkpq}\\
&+\Big[-\frac{4s(n^2-3n+4)+4(n-2)}{n(n-2)^2}|E|^2+\frac{2s}{n}|W|^2\Big]g_{ij}\\
&+\frac{4s(n^2-3n+4)+4(n-2)}{(n-2)^2}E_{ik}E_{jk}\\
&+\frac{4-2n-2n(n-1)t+4(n-2)s}{n(n-1)}RE_{ij},
\endaligned\end{equation}
and the desired estimate \eqref{2lem-form1} follows from
\begin{equation}\label{add2-Sec-3}\aligned
\frac{1+4s}{2}\Delta |E|^2=&(1+4s)|\nabla E|^2+(1+4s)E_{ij}\Delta E_{ij}.
\endaligned\end{equation}

By virtue of the Ricci identity and the formula \eqref{2-Sec-15}, we have
\begin{equation}\label{2-Sec-13}\aligned
E_{kj,ik}=&E_{kj,ki}+E_{lj}R_{lkik}+E_{kl}R_{ljik}\\
=&E_{kj,ki}+E_{lj}R_{li}+E_{kl}R_{ljik}\\
=&\frac{n-2}{2n}R_{,ij}+E_{ik}E_{jk}+\frac{1}{n}RE_{ij}\\
&-\Big[E_{kl}W_{ikjl}+\frac{1}{n-2}(|E|^2g_{ij}-2E_{ik}E_{jk})-\frac{1}{n(n-1)}R E_{ij}\Big]\\
=&\frac{n-2}{2n}R_{,ij}+\frac{n}{n-2}E_{ik}E_{jk}+\frac{1}{n-1}RE_{ij}\\
&-E_{kl}W_{ikjl}-\frac{1}{n-2}|E|^2g_{ij},
\endaligned\end{equation}
where we used the second Bianchi identity $E_{kj,k}=\frac{n-2}{2n}R_{,j}$.
Thus, from \eqref{2-Sec-5}, we have
\begin{equation}\label{2-Sec-14}\aligned
(1+4s)C_{kij,k}=&(1+4s)\Delta E_{ij}-(1+4s)E_{kj,ik}+\frac{(n-2)(1+4s)}{2n(n-1)}[(\Delta R)g_{ij}-R_{,ij}]\\
=&(1+4s)\Delta E_{ij}-(1+4s)\Big[\frac{n-2}{2n}R_{,ij}+\frac{n}{n-2}E_{ik}E_{jk}+\frac{1}{n-1}RE_{ij}\\
&-E_{kl}W_{ikjl}-\frac{1}{n-2}|E|^2g_{ij}\Big]+\frac{(n-2)(1+4s)}{2n(n-1)}[(\Delta R)g_{ij}-R_{,ij}]\\
=&(1+4s)\Delta E_{ij}-\frac{(n-2)(1+4s)}{2(n-1)}\mathring{R}_{,ij}-(1+4s)\Big[\frac{n}{n-2}E_{ik}E_{jk}\\
&+\frac{1}{n-1}RE_{ij}-E_{kl}W_{ikjl}-\frac{1}{n-2}|E|^2g_{ij}\Big].
\endaligned\end{equation}
Inserting \eqref{add2-Sec-2} into the inequality \eqref{2-Sec-14} yields
the estimate \eqref{2lem-form2}.

\section{Proof of main theorems}
\subsection{Proof of Theorem \ref{thm1-1}. }

Under the condition that $n\neq 4$, we have from \eqref{2lem-form2}
\begin{equation}\label{2-Sec-19}\aligned
0=&\int_M\Big[-(1+4s)C_{kij,k}E_{ij}+\frac{n+4(n-1)t+4s}{2(n-1)}\mathring{R}_{,ij}E_{ij}\\
&-\frac{n-2+8s}{n-2}W_{ikjl}E_{kl}E_{ij}-2sW_{ikpq}W_{jkpq}E_{ij}\\
&-\frac{(n-4)(n-2+4s)}{(n-2)^2}E_{ij}E_{jk}E_{ki}+\frac{4-3n-8s-2n(n-1)t}{n(n-1)}R |E|^2\Big]\,dv.
\endaligned\end{equation}
Since
$$\int_M \mathring{R}_{,ij}E_{ij}\,dv=\int_M R_{,ij}E_{ij}\,dv=-\int_M R_{,i}E_{ij,j}\,dv=-\frac{n-2}{2n}\int_M |\nabla R|^2 \,dv$$
and
$$\int_MC_{kij,k}E_{ij}\,dv=-\int_MC_{kij}E_{ij,k}\,dv=-\frac{1}{2}\int_M|C|^2\,dv,$$
then \eqref{2-Sec-19} can be written as
\begin{equation}\label{2-Sec-20}\aligned
0=&\int_M\Big[\frac{1+4s}{2}|C|^2-\frac{(n-2)[n+4(n-1)t+4s]}{4n(n-1)}|\nabla R|^2\\
&-\frac{(n-4)(n-2+4s)}{(n-2)^2}E_{ij}E_{jk}E_{ki}+\frac{4-3n-8s-2n(n-1)t}{n(n-1)}R |E|^2\\
&-\frac{n-2+8s}{n-2}W_{ikjl}E_{kl}E_{ij}-2sW_{ikpq}W_{jkpq}E_{ij}\Big]\,dv.
\endaligned\end{equation}
In particular, noticing that
\begin{equation}\label{addProof11}
-\frac{n-2}{\sqrt{n(n-1)}}|E|^3\leq E_{ij}E_{jk}E_{ki}\leq\frac{n-2}{\sqrt{n(n-1)}}|E|^3
\end{equation}
with equality in \eqref{addProof11} at some point $p\in M$ if and
only if $E$ can be diagonalized at $p$ and the eigenvalue
multiplicity of $E$ is at least $n-1$. If $|E|\neq0$ and the equality in the left hand side of  \eqref{addProof11} occurs, then $n-1$ of eigenvalues which are equal must be positive (see \cite{Okumura} or Lemma 5.1 in \cite{HL2004}).
On the other hand, the following inequality which was first proved by Huisken
(cf. \cite[Lemma 3.4]{Huisken1985}):
\begin{equation}\label{addProof10}
|W_{ikjl}E_{ij}E_{kl}|\leq\sqrt{\frac{n-2}{2(n-1)}}|W||E|^2.
\end{equation}
When $n+4(n-1)t+4s>0$, applying \eqref{addProof10} and \eqref{addProof11} into \eqref{2-Sec-20} yields
\begin{equation}\label{2-Sec-22}\aligned
0\leq&\int_M\Big[\frac{1+4s}{2}|C|^2+\frac{(n-4)|n-2+4s|}{(n-2)\sqrt{n(n-1)}}|E|^3\\
&+\frac{4-3n-8s-2n(n-1)t}{n(n-1)}R |E|^2\\
&+\frac{|n-2+8s|}{\sqrt{2(n-2)(n-1)}}|W||E|^2+2|s||W|^2|E|\Big]\,dv.
\endaligned\end{equation}
If $n+4(n-1)t+4s<0$, then from \eqref{2-Sec-20}, we have
\begin{equation}\label{2-Sec-23}\aligned
0\geq&\int_M\Big[\frac{1+4s}{2}|C|^2-\frac{(n-4)|n-2+4s|}{(n-2)\sqrt{n(n-1)}}|E|^3\\
&+\frac{4-3n-8s-2n(n-1)t}{n(n-1)}R |E|^2\\
&-\frac{|n-2+8s|}{\sqrt{2(n-2)(n-1)}}|W||E|^2-2|s||W|^2|E|\Big]\,dv.
\endaligned\end{equation}
In particular, equalities in \eqref{2-Sec-22} or \eqref{2-Sec-23} occur if and only if $R$ is constant.
Hence, as stated in the lines following
\eqref{addProof11}, $E$ has, at each point $p$, an eigenvalue of
multiplicity $n-1$ or $n$. For $n=3$, it is well known that $W=0$
and \eqref{addProof10} is an equality. When $n\geq 5$, writing
$E_{ij}=ag_{ij}+bv_iv_j$ at $p$, with two scalars $a,b$ and a vector
$v$, we see that the left-hand side of \eqref{addProof10} is zero at
$p$.  This shows that $M^n$, $n\geq5$, must
be conformally flat or Einstein due to the equality in \eqref{addProof10}.

Furthermore, for the case of $n+4(n-1)t+4s>0$, if the equality in
\eqref{2-Sec-22} occurs and $M^n$ is not Einstein (that is,
$E\neq0$), we have $W=0$ according to above arguments, which shows
from \eqref{2-Sec-22} that
\begin{equation}\label{addadd2Propformula2}\aligned
\int_M\Big[\frac{1+4s}{2}|C|^2&+\frac{(n-4)|n-2+4s|}{(n-2)\sqrt{n(n-1)}}|E|^3\\
&+\frac{4-3n-8s-2n(n-1)t}{n(n-1)}R |E|^2\Big]\,dv=0.
\endaligned\end{equation}
In particular, from \eqref{2-Sec-3}, we have that $|E|^2$ is constant if $s\neq-\frac{n-2}{4}$. As a result, the eigenvalues of Ricci curvature are constant from that both $E$ and $R$ are constant, which shows $C_{ijk}=0$. In this case, $\nabla E=0$ and from the de Rham
decomposition theorem, $M^n$ splits as a product of two Einstein manifolds
$N^1\times N^{n-1}$, where $N^{n-1}$ is a Einstein manifold. Let $\lambda_1,\cdots,\lambda_n$ be the eigenvalues of Ricci curvature with $\lambda_2=\cdots=\lambda_n$. Since the dimension of $N^1$ is one, we have $\lambda_1=0$ and $R=\sqrt{n(n-1)}|E|$. Thus, \eqref{addadd2Propformula2} becomes
\begin{equation}\label{March-Prf2-2}\aligned
\Big[\frac{(n-4)|n-2+4s|}{(n-2)}+4-3n-8s-2n(n-1)t\Big]\int_M\frac{|E|^3}{\sqrt{n(n-1)}}\,dv=0,
\endaligned\end{equation}
which shows $\frac{(n-4)|n-2+4s|}{(n-2)}+4-3n-8s-2n(n-1)t=0$.

The proof of the case (2) in Theorem \ref{thm1-1} is similar, we omit it here.
It completes the proof of Theorem \ref{thm1-1}.\endproof

\subsection{Proof of Theorem \ref{thm1-2}. }
We recall that if $E$ is a Codazzi tensor, then it
satisfies the following sharp inequality (for the proof, for
instance, see \cite{Hebey1996}. This inequality was first observed by
Bourguignon \cite{Bourguignon1990}):
\begin{equation}\label{3Eq-1}
|\nabla E|^2\geq \frac{n+2}{n}|\nabla|E||^2,
\end{equation} and $R$ is constant.
Inserting \eqref{3Eq-1}, \eqref{addProof11} and  \eqref{addProof10} into \eqref{2lem-form1} yields
\begin{equation}\label{3Eq-2}\aligned
\frac{1+4s}{2}\Delta |E|^2\geq&\frac{(n+2)(1+4s)}{n}|\nabla|E||^2
-\frac{|2(n-2)+4ns|}{\sqrt{2(n-2)(n-1)}}|W||E|^2\\
&-2|s||W|^2|E|-\frac{|4s(n^2-3n+4)+4(n-2)|}{(n-2)\sqrt{n(n-1)}}|E|^3\\
&+\frac{4-2n-2n(n-1)t+4(n-2)s}{n(n-1)}R|E|^2
\endaligned\end{equation} from $1+4s\geq0$.

Let $\Omega_0=\{p\in M:\,|E|\neq0\}$.  By virtue of the Lemma 2.2 in \cite{Catino2016} (or
see \cite[Theorem 1.8]{Kazdan1998}), one has ${\rm Vol}(M\backslash
\Omega_0)=0$. For any $\epsilon>0$, we define
$\Omega_{\epsilon}=\{p\in M:|E|\geq\epsilon\}$ and
$$f_{\epsilon}(p)=\left\{\begin{array}{ll}
|E|(p)& \ {\rm if} \ \ p\in \Omega_{\epsilon};\\
\varepsilon& \ {\rm if} \ \ p\in M\backslash\Omega_{\epsilon}.
\end{array}\right.$$
Then at the regular value $\epsilon$ of $|E|$, we have
parts
\begin{equation}\label{3Eq-3}\aligned
\int_{M}&\Bigg(-\frac{1}{2}\Delta
|E|^2+\frac{n+2}{n}|\nabla|E||^2\Bigg)f_{\epsilon}^{-\frac{n+2}{n}}\\
=&-\frac{n+2}{n}\int_{M}\langle\nabla |E|,\nabla
f_{\epsilon}\rangle|E|f_{\epsilon}^{-\frac{n+2}{n}-1}
+\frac{n+2}{n}\int_{M}|\nabla|E||^2f_{\epsilon}^{-\frac{n+2}{n}}\\
\endaligned\end{equation} which tends to the zero as $\epsilon \rightarrow 0$, where in the last equality we used
$f_{\epsilon}=|E|$ on $\Omega_{\epsilon}$ and $\nabla
f_{\epsilon}=0$ on $M\backslash \Omega_{\epsilon}$.
Multiplying both sides of inequality \eqref{3Eq-2} by
$f_{\epsilon}^{-\frac{n+2}{n}}$ and noticing
\eqref{3Eq-3}, we derive
\begin{equation}\label{3Eq-4}\aligned
0\geq&\int_{M}\Big[-\frac{|2(n-2)+4ns|}{\sqrt{2(n-2)(n-1)}}|W||E|^2-2|s||W|^2|E|\\
&-\frac{|4s(n^2-3n+4)+4(n-2)|}{(n-2)\sqrt{n(n-1)}}|E|^3\\
&+\frac{4-2n-2n(n-1)t+4(n-2)s}{n(n-1)}R|E|^2\Big]f_{\epsilon}^{-\frac{n+2}{n}}\\
=&\int_{M}\Big[-\frac{|2(n-2)+4ns|}{\sqrt{2(n-2)(n-1)}}|W||E|^{\frac{n-2}{n}}-2|s||W|^2|E|^{-\frac{2}{n}}\\
&-\frac{|4s(n^2-3n+4)+4(n-2)|}{(n-2)\sqrt{n(n-1)}}|E|^{\frac{2(n-1)}{n}}\\
&+\frac{4-2n-2n(n-1)t+4(n-2)s}{n(n-1)}R|E|^{\frac{n-2}{n}}\Big]|E|^{\frac{n+2}{n}}f_{\epsilon}^{-\frac{n+2}{n}}.
\endaligned\end{equation}
Letting $\epsilon\rightarrow 0$ in \eqref{3Eq-4}, we obtain
$|E|^{\frac{n+2}{n}}f_{\epsilon}^{-\frac{n+2}{n}}\rightarrow 1$
a.e. on $M^n$ and the desired estimate \eqref{3For-th-1}
follows. Moreover, we can obtain the case of equality in \eqref{3For-th-1} by using the arguments as in Theorem \ref{thm1-1}.

\subsection{Proof of Theorem \ref{thm1-3}. }

In order to prove Theorem \ref{thm1-3}, we fist prove the following result:

\begin{lem}\label{3lem1}
Let $M^{n}$ be a locally conformally flat closed manifold of dimension $n\neq4$ with
positive scalar curvature and $g$ be a critical metric for
$\mathcal{F}_{t,s}$ on $\mathscr{M}_1(M^{n})$. If for $n=3$, $t,s$ satisfy
\begin{equation}\label{3lem-Form-1}
\begin{cases}1+2t+2s>0\\
(1+4s)(3+8t+4s)>0;
\end{cases}
\end{equation}
for $n\geq5$, $t,s$ satisfy
\begin{equation}\label{3lem-Form-1-1}
\begin{cases}1+2t+2s>0\\
[(n-2)+4s][n+4(n-1)t+4s]<0,
\end{cases}
\end{equation}
then we have
\begin{equation}\label{3lem-Form-2}\aligned
\int_M&\Big[(1+4s)R^2+\Big(\frac{1+4s}{2}+\frac{(n-1)(n-4)(n-2+4s)(1+2t+2s)}{n(n-2)[n+4(n-1)t+4s]}\Big)^2\\
&\times\frac{2n^2[n+4(n-1)t+4s]}{(n-2)(n-4)(1+2t+2s)[(n-2)+4s]}|E|^2\Big]\frac{|\nabla E|^2}{R}\,dv\\
\leq&\int_M\Big[\frac{|4s(n^2-3n+4)+4(n-2)|}{(n-2)\sqrt{n(n-1)}}|E|\\
&-\frac{4-2n-2n(n-1)t+4(n-2)s}{n(n-1)}R\Big]R|E|^2\,dv.
\endaligned\end{equation}
\end{lem}

\proof From the well-known Bochner formula, we have
$$\aligned
\frac{1}{2}\Delta|\nabla R|^{2}=&|\nabla^{2} R|^{2}+{\rm Ric}(\nabla R, \nabla R)+\langle\nabla\Delta R,\nabla R \rangle\\
=&|\nabla^{2} R|^{2}+E(\nabla R, \nabla R)+\frac{1}{n}R |\nabla
R|^{2}+\langle \nabla\Delta R, \nabla R\rangle\\
\geq&\frac{1}{n}(\Delta R)^{2}+E(\nabla R, \nabla R)+\frac{1}{n}R |\nabla
R|^{2}+\langle \nabla\Delta R, \nabla R\rangle,
\endaligned$$ where the Cauchy
inequality $|\nabla^{2}R|^{2}\geq\frac{1}{n}(\Delta R)^2$ is used.
Integrating both sides of the above inequality yields
\begin{equation}\label{3lem-Form-3}
\int_ME(\nabla R, \nabla
R)\,dv\leq\frac{n-1}{n}\int_M(\Delta
R)^2\,dv-\frac{1}{n}\int_MR |\nabla R|^{2}\,dv.
\end{equation} Using the equation \eqref{2-Sec-3} with $n+4(n-1)t+4s\neq0$, one has
\begin{equation}\label{3lem-Form-4}\aligned
\int_M(\Delta R)^{2}\,dv=&\int_M R \Delta^{2} R \,dv\\
=&\frac{n-4}{n+4(n-1)t+4s}\int_M R \Big(\frac{n-2+4s}{n-2}\Delta |E|^{2}\,dv\\
&+\frac{n-1+n(n-1)t+2s}{n(n-1)}\Delta R^{2}\Big)\,dv\\
=&-\frac{n-4}{n+4(n-1)t+4s}\int_M\Big(\frac{n-2+4s}{n-2}\langle \nabla|E|^{2},
\nabla R\rangle\\
&+\frac{2[n-1+n(n-1)t+2s]}{n(n-1)}R|\nabla R|^{2}\Big)\,dv.
\endaligned\end{equation}
Putting \eqref{3lem-Form-4} into \eqref{3lem-Form-3} yields
\begin{equation}\label{3lem-Form-5}\aligned
\int_ME(\nabla R, \nabla R)\,dv
\leq&-\frac{(n-1)(n-4)(n-2+4s)}{n(n-2)[n+4(n-1)t+4s]}\int_M\langle
\nabla|E|^{2}, \nabla R\rangle\,dv\\
&-\Big[\frac{2(n-4)[n-1+n(n-1)t+2s]}{n^2[n+4(n-1)t+4s]}+\frac{1}{n}\Big]\int_MR
|\nabla R|^{2}\,dv.
\endaligned\end{equation}
Since $W=0$, multiplying both sides of \eqref{2lem-form1} with $R$ and integrating it, we have
\begin{equation}\label{3lem-Form-6}\aligned
\int_M&\Big[\frac{1+4s}{2}\langle\nabla|E|^{2}, \nabla R\rangle+(1+4s)R|\nabla E|^2-\frac{(n-2)(1+2t+2s)}{2n}R|\nabla R|^2\\
&-(1+2t+2s)E(\nabla R, \nabla R)\Big]\,dv\\
=&\int_M\Big[-\frac{4s(n^2-3n+4)+4(n-2)}{(n-2)^2}RE_{ik}E_{kj}E_{ji}\\
&-\frac{4-2n-2n(n-1)t+4(n-2)s}{n(n-1)}R^2|E|^2\Big]\,dv\\
\leq&\int_M\Big[\frac{|4s(n^2-3n+4)+4(n-2)|}{(n-2)\sqrt{n(n-1)}}|E|\\
&-\frac{4-2n-2n(n-1)t+4(n-2)s}{n(n-1)}R\Big]R|E|^2\,dv.
\endaligned\end{equation}
Notcing $1+2t+2s>0$, inserting \eqref{3lem-Form-5} into \eqref{3lem-Form-6} gives
\begin{equation}\label{3lem-Form-7}\aligned
\int_M&\Big[\Big(\frac{1+4s}{2}+\frac{(n-1)(n-4)(n-2+4s)(1+2t+2s)}{n(n-2)[n+4(n-1)t+4s]}\Big)\langle\nabla|E|^{2}, \nabla R\rangle\\
&+(1+2t+2s)\Big(\frac{2(n-4)[n-1+n(n-1)t+2s]}{n^2[n+4(n-1)t+4s]}+\frac{1}{n}\\
&-\frac{n-2}{2n}\Big)R|\nabla R|^2+(1+4s)R|\nabla E|^2\Big]\,dv\\
\leq&\int_M\Big[\frac{|4s(n^2-3n+4)+4(n-2)|}{(n-2)\sqrt{n(n-1)}}|E|\\
&-\frac{4-2n-2n(n-1)t+4(n-2)s}{n(n-1)}R\Big]R|E|^2\,dv.
\endaligned\end{equation}
For any positive constant $\varepsilon$, it holds that
\begin{equation}\label{3lem-Form-8}\aligned
&\Big(\frac{1+4s}{2}+\frac{(n-1)(n-4)(n-2+4s)(1+2t+2s)}{n(n-2)[n+4(n-1)t+4s]}\Big)\langle\nabla|E|^{2}, \nabla R\rangle\\
\geq&-2\Big|\frac{1+4s}{2}+\frac{(n-1)(n-4)(n-2+4s)(1+2t+2s)}{n(n-2)[n+4(n-1)t+4s]}\Big||E||\nabla
|E|||\nabla R|\\
\geq&-2\Big|\frac{1+4s}{2}+\frac{(n-1)(n-4)(n-2+4s)(1+2t+2s)}{n(n-2)[n+4(n-1)t+4s]}\Big||E||\nabla
E||\nabla R|\\
\geq&-\Big|\frac{1+4s}{2}+\frac{(n-1)(n-4)(n-2+4s)(1+2t+2s)}{n(n-2)[n+4(n-1)t+4s]}\Big|\\
&\times\Big[\varepsilon
R|\nabla R|^2+\frac{1}{\varepsilon}\frac{|E|^2}{R}|\nabla E|^2\Big],
\endaligned\end{equation}
where, in the third line, we use the Kato inequality $|\nabla E|\geq|\nabla|E||$. Hence, when $t,s$ satisfy \eqref{3lem-Form-1} or \eqref{3lem-Form-1-1}, we both have
$$\aligned
(1+2t+2s)&\Big(\frac{2(n-4)[n-1+n(n-1)t+2s]}{n^2[n+4(n-1)t+4s]}+\frac{1}{n}-\frac{n-2}{2n}\Big)\\
=&-\frac{(n-2)(n-4)(1+2t+2s)[(n-2)+4s]}{2n^2[n+4(n-1)t+4s]}\\
>&0.
\endaligned$$
Therefore, there exists a positive constant $\varepsilon_0$ such that
\begin{equation}\label{3lem-Form-9}\aligned
(1+&2t+2s)\Big(\frac{2(n-4)[n-1+n(n-1)t+2s]}{n^2[n+4(n-1)t+4s]}+\frac{1}{n}-\frac{n-2}{2n}\Big)\\
-&\Big|\frac{1+4s}{2}+\frac{(n-1)(n-4)(n-2+4s)(1+2t+2s)}{n(n-2)[n+4(n-1)t+4s]}\Big|\varepsilon_0=0.
\endaligned\end{equation}
Inserting \eqref{3lem-Form-8} with $\varepsilon_0$ into \eqref{3lem-Form-7}
yields the desired estimate \eqref{3lem-Form-2} and the proof of Lemma \ref{3lem1} is completed.\endproof

Now, we will complete the proof of Theorem \ref{thm1-3} with the help of the Lemma \ref{3lem1}. When
$n=3$, \eqref{3lem-Form-2} becomes
\begin{equation}\label{3lem-Form-10}\aligned
(1+4s)\int_M&\Big[R^2-
\frac{(5+16t+4s)^2}{2(3+8t+4s)(1+2t+2s)}|E|^2\Big]\frac{|\nabla E|^2}{R}\,dv\\
\leq&\int_M\Big[\frac{4}{\sqrt{6}}|1+4s||E|+\frac{1+6t-2s}{3}R\Big]R|E|^2\,dv.
\endaligned\end{equation}
Thus, under the assumption \eqref{3lem-Form-11} and \eqref{3lem-Form-12}, we have
\begin{equation}\label{3lem-Form-15}\aligned
0\leq&(1+4s)\int_M\Big[R^2-
\frac{(5+16t+4s)^2}{2(3+8t+4s)(1+2t+2s)}|E|^2\Big]\frac{|\nabla E|^2}{R}\,dv\\
\leq&\int_M\Big[\frac{4}{\sqrt{6}}|1+4s||E|+\frac{1+6t-2s}{3}R\Big]R|E|^2\,dv\\
\leq&0,
\endaligned\end{equation}
which gives $E=0$ and $M^3$ is Einstein. Thus, $M^3$ is of constant positive sectional curvature.

When $n\geq5$, under the assumption \eqref{3lem-Form-21} and \eqref{3lem-Form-22},
we have
\begin{equation}\label{3lem-Form-25}\aligned
0\leq&\int_M\Big[(1+4s)R^2+\Big(\frac{1+4s}{2}+\frac{(n-1)(n-4)(n-2+4s)(1+2t+2s)}{n(n-2)[n+4(n-1)t+4s]}\Big)^2\\
&\times\frac{2n^2[n+4(n-1)t+4s]}{(n-2)(n-4)(1+2t+2s)[(n-2)+4s]}|E|^2\Big]\frac{|\nabla E|^2}{R}\,dv\\
\leq&\int_M\Big[\frac{|4s(n^2-3n+4)+4(n-2)|}{(n-2)\sqrt{n(n-1)}}|E|\\
&-\frac{4-2n-2n(n-1)t+4(n-2)s}{n(n-1)}R\Big]R|E|^2\,dv\\
\leq&0,
\endaligned\end{equation}
which gives $E=0$ and $M^n$ is Einstein. Thus, $M^n$ is also of constant positive sectional curvature.


\bibliographystyle{Plain}

\begin{thebibliography}{10}

\bibitem{Anderson1997}
M.  Anderson, Extrema of curvature functionals on the space of
metrics on 3-manifolds, Calc. Var. Partial Differential Equations 5
(1997), 199-269.

\bibitem{Barros2017}
A. Barros, A. Da Silva, On locally conformally flat critical metrics for quadratic
functionals, Ann. Glob. Anal. Geom. 52 (2017), 1-9.

\bibitem{Bess2008} A. Besse, Einstein Manifolds, Springer-Verlag,
Berlin, 2008.

\bibitem{Bourguignon1990}
J.-P. Bourguignon, The ``magic" of Weitzenb\"{o}ck formulas,
Variational methods (Paris, 1988), 251-271, Progr. Nonlinear
Differential Equations Appl., 4, Birkh\"{a}user Boston, Boston, MA,
1990.

\bibitem{Catino2016}
G. Catino, On conformally flat manifolds with constant positive
scalar curvature, Proc. Amer. Math. Soc. 144 (2016), 2627-2634.

\bibitem{Catino2015} G. Catino, Some rigidity results on critical metrics for quadratic
functionals,  Calc. Var. Partial Differential Equations, 54 (2015),
 2921-2937.

\bibitem{Rham1952}
G. de Rham, Sur la reductibilit\'{e} d'un espace de Riemann,
Comment. Math. Helv. 26, (1952). 328-344.

\bibitem{Gursky2001}
M. Gursky, J. Viaclovsky, A new variational characterization of
three-dimensional space forms, Invent. Math. 145 (2001), 251-278.

\bibitem{Gursky2015}
M. Gursky, J. Viaclovsky, Rigidity and stability of Einstein metrics
for quadratic curvature functionals, J. Reine Angew. Math. 700
(2015), 37-91.

\bibitem{Hebey1996}
E. Hebey, M. Vaugon, Effective $L_p$ pinching for the concircular
curvature, J. Geom. Anal. 6 (1996), 531-553.

\bibitem{HC2017}
G. Huang, L. Chen, Some characterizations on critical metrics for quadratic curvature
functions, to appear in {\it  Proc. Amer. Math. Soc.}

\bibitem{Huisken1985}
G. Huisken, Ricci deformation of the metric on a Riemannian
manifold. J. Differ. Geom. 21 (1985), 47-62.

\bibitem{HL2004}
Z. Hu,  H. Li, A new variational characterization of
$n$-dimensional space forms, Trans. Amer. Math. Soc. 356 (2004),
3005-3023.

\bibitem{Hu2010}
Z. Hu, S. Nishikawa, U. Simon, Critical metrics of the Schouten
functional, J. Geom. 98 (2010), 91-113.

\bibitem{Kazdan1998}
J. Kazdan, Unique continuation in geometry, Comm. Pure Apple. Math.
41 (1988), 667-681.



\bibitem{Lamontagne1994}
F. Lamontagne, Une remarque sur la norme $L^2$ du tenseur de
courbure, C. R. Acad. Sci. Paris S\'{e}r. I Math. 319 (1994),
237-240.

\bibitem{Lamontagne1998}
F. Lamontagne,  A critical metric for the $L^2$-norm of the
curvature tensor on $\mathbb{S}^3$, Proc. Amer. Math. Soc. 126
(1998), 589-593.

\bibitem{Okumura} M. Okumura, Hypersurfaces and a pinching problem
on the second fundamental tensor, Amer. J. Math. 96 (1974), 207-213.

\end{thebibliography}

\end{document}